\documentclass[fleqn,12pt]{article}
\usepackage{amsfonts}
\usepackage{amsthm}
\usepackage{amssymb}
\usepackage{amsmath}

\newtheorem*{theo}{\bf Theorem}
\newtheorem{lemma}{\bf Lemma}
\newtheorem{coro}{\bf Corollary}

\newcommand{\Hex}{{\cal H}}
\newcommand{\Z}{{\mathbb Z}}

\newcommand{\R}{{\mathbb R}}

\newcommand{\bea}{\begin{eqnarray*}}
\newcommand{\eea}{\end{eqnarray*}}
\newcommand{\be}{\begin{eqnarray}}
\newcommand{\ee}{\end{eqnarray}}

\newcommand{\vol}{\mbox{vol}\,}

\newcommand{\ve}{\bf}

\begin{document}

\setlength{\unitlength}{1mm}
\setcounter{page}{1}

\title{Siegel's Lemma and Sum--Distinct Sets}
\author{by Iskander ALIEV}
\date{}
\maketitle \vskip 1cm

\begin{center}

School of Mathematics

University of Edinburgh

James Clerk Maxwell Building

King's Buildings, Mayfield Road

Edinburgh EH9 3JZ, UK

I.Aliev@ed.ac.uk

Tel.: +44 131 650 5056

Fax: +44 131 650 6553

\end{center}

\newpage

\noindent {\bf Abstract.}
Let $L({\ve x}) = a_1x_1 + a_2x_2 + \ldots + a_nx_n$, $n\ge 2$ be a
linear form with integer coefficients $a_1, a_2, \ldots , a_n$ which
are not all zero. A basic problem is to determine nonzero integer
vectors ${\ve x}$ such that $L({\ve x}) = 0$, and the maximum norm
$||{\ve x}||$ is relatively small compared with the size of the
coefficients $a_1, a_2, \ldots , a_n$.
%
%
The main result of the paper asserts that 
there exist linearly independent vectors ${\ve x}_1,\ldots, {\ve
x}_{n-1}\in\Z^n$ such that $L({\ve x}_i)=0$, $i=1,\ldots, n-1$ and
$$ ||{\ve x}_1||\cdots||{\ve x}_{n-1}||<\frac{||{\ve
a}||}{\sigma_n}\,,$$
%
%
where ${\ve a}=(a_1, a_2, \ldots , a_n)$ and
$$\sigma_n=\frac{2}{\pi}\int_0^\infty\left(\frac{\sin t}{t}\right)^n
dt\,. $$

This result also implies a new lower bound on the greatest element
of a sum--distinct set of positive integers (Erd\"os--Moser
problem).
The main tools are the Minkowski theorem on successive minima and
the Busemann theorem from convex geometry.

\bigskip
\noindent {\bf Keywords:} sections of the cube, sinc integrals,
Busemann's theorem, intersection body, successive minima

\bigskip
\noindent {\bf 2000 MS Classification:} 11H06, 11P70

\newpage

\section{Introduction}

Let ${\ve a}=(a_1,\ldots,a_n)$, $n\ge 2$ be a non--zero integral
vector. Consider the linear form $L({\ve x}) = a_1x_1 + a_2x_2 +
\ldots + a_nx_n$.
%
%
Siegel's Lemma w. r. t. the maximum norm $||\cdot||$ asks for an
optimal constant $c_n>0$ such that the equation
%
\bea L({\ve x})=0\, \label{ax}\eea
has an integral solution ${\ve x}=(x_1,\ldots,x_n)$ with
\be 0<||{\ve x}||^{n-1}\le c_n||{\ve a}||\,. \label{c_n}\ee
%

%
%
The only known exact values of $c_n$ are $c_2=1$, $c_3=4/3$ and
$c_4=27/19$ (see \cite{Disser}, \cite{NewSL}). Note that for $n=3,4$
the equality in (\ref{c_n}) is not attained. A. Schinzel
\cite{NewSL} has shown that for $n\ge 3$
\bea c_n=\sup \Delta( \Hex^{n-1}_{\alpha_1,\ldots,\alpha_{n-3}}
)^{-1}\ge 1\,, \eea where $\Delta(\cdot)$ denotes the critical
determinant, $\Hex^{n-1}_{\alpha_1,\ldots,\alpha_{n-3}}$ is a
generalized hexagon in $\R^{n-1}$ given by \bea |x_i|\le
1\,,\;\;\;i=1,\ldots,n-1\,,\;\;\;
|\sum_{i=1}^{n-3}\alpha_ix_i+x_{n-2}+x_{n-1}|\le 1 \eea and
$\alpha_i$ range over all rational numbers in the interval
$(\,0\,,\;1\,]$.
The values of $c_n$ for $n\le 4$ indicate that, most likely,
$c_n=\Delta( \Hex^{n-1}_{1,\ldots,1} )^{-1}$. However, a proof of
this conjecture does not seem within reach at present.
The best known upper bound
\be c_n\le\sqrt{n}\,\label{old_bound}\ee
follows from the classical result of Bombieri and Vaaler
(\cite{BombVaal}, Theorem 1).

In the present paper we estimate $c_n$ via values of the {\em
sinc} integrals
\bea \sigma_n=\frac{2}{\pi}\int_0^\infty\left(\frac{\sin
t}{t}\right)^n dt\,. \eea
The main result is as follows:
\begin{theo}
For any non--zero vector ${\ve a}\in\Z^n$, $n\ge 5$, there exist
linearly independent vectors ${\ve x}_1,\ldots, {\ve
x}_{n-1}\in\Z^n$ such that $L({\ve x}_i)=0$, $i=1,\ldots, n-1$ and
\be ||{\ve x}_1||\cdots||{\ve x}_{n-1}||<\frac{||{\ve
a}||}{\sigma_n}\,.\label{product}\ee \label{main}
\end{theo}
From (\ref{product}) we immediately get the bound \be
c_n\le\sigma_n^{-1}\,,\label{tozhe_nasha}\ee and since
\be \sigma_n^{-1}\sim\sqrt{\frac{\pi n}{6}}\,,\;\;\;\mbox{as}\,\;
n\rightarrow\infty\,\label{asympt_sigma}\ee
(see Section \ref{sinc}), the theorem 
asymptotically
improves the estimate (\ref{old_bound}). It is also known (see e. g.
\cite{MeRo}) that
\bea \sigma_n=\frac{n}{2^{n-1}}\sum\limits_{0\le
r<n/2,\,r\in\Z}\frac{(-1)^r(n-2r)^{n-1}}{r!(n-r)!}\,.\eea
The sequences of numerators and denominators of $\sigma_n/2$ can be
found in \cite{Sloane}.
%


\noindent{{\bf Remark 1}}

\begin{enumerate}
\item[(1)] Calculation shows that for all $5\le n\le 1000$
the bound (\ref{tozhe_nasha}) is slightly better than
(\ref{old_bound}).

\item[(2)] For $n\le 4$ the constant $\sigma_n^{-1}$ in (\ref{product})
can be replaced by $c_n$.
This follows from the observation that any origin--symmetric convex
body in $\R^n$, $n\le 3$ has anomaly $1$ (see \cite{Woods}).

\end{enumerate}

As it was observed by A. Schinzel (personal communication), Siegel's
Lemma w. r. t. maximum norm can be applied to the following well
known problem from additive number theory. A finite set
$\{a_1,\ldots,a_n\}$ of integers is called {\em sum--distinct set}
if any two of its $2^n$ subsums differ by at least $1$. We shall
assume w. l. o. g. that $0<a_1<a_2<\ldots<a_n$. In 1955, P. Erd\"os
and L. Moser (\cite{EM}, Problem 6) asked for an estimate on the
least possible $a_n$ of such a set. They proved that
\be a_n>\max\left \{\frac{2^n}{n}, \frac{2^n}{4\sqrt{n}}\right
\}\label{EM_result}\ee
and Erd\"os conjectured that $a_n>C_0 2^n$, $C_0>0$. In 1986, N. D.
Elkies \cite{El} showed that
\be a_n>2^{-n} {2n\choose n}\label{Elk}\ee
and this result is still cited by Guy (\cite{Guy}, Problem C8) as
the best known lower bound for large $n$.
Following \cite{El}, note that references \cite{EM, Guy} state the
problem equivalently in terms of ,,inverse function''. They ask to
maximize the size $m$ of a sum--distinct subset of $\{1,2,\ldots,
x\}$, given $x$. Clearly, the bound $a_n>C_1 n^{-s}2^n$ corresponds
to
\bea m<\log_2 x+s \log_2\log_2 x+\log_2\frac{1}{C_1}-o(1)\,. \eea

\begin{coro}
For any sum--distinct set $\{a_1,\ldots,a_n\}$ with
$0<a_1<\ldots<a_n$, the inequality
\be a_n>\sigma_n 2^{n-1} \label{nasha}\ee
holds.
\label{lower}
\end{coro}
%
%
Since 
\bea 2^{-n}{2n\choose n}\sim \frac{2^n}{\sqrt{\pi
n}}\, \;\;\;
\mbox{and}\,\;\;\;
\sigma_n 2^{n-1}\sim \frac{2^n}{\sqrt{\frac{2\pi n}{3}}}\,,
\;\;\;\mbox{as}\; n\rightarrow \infty\,,
\eea
Corollary \ref{lower} asymptotically improves the result of Elkies
with factor $\sqrt{3/2}$.

\noindent{{\bf Remark 2}}

\begin{enumerate}

\item[(1)]
Sum--distinct sets with minimal largest element are known up to
$n=9$ (see \cite{BoMo}). In the latter case the estimate
(\ref{nasha}) predicts $a_9\ge 116$ and the optimal bound is $a_9\ge
161$.
Calculation shows that for all $10\le n\le 1000$ the bound
(\ref{nasha}) is slightly  better than (\ref{Elk}).

\item[(2)] Prof. Noam Elkies kindly informed the author about existing of
an unpublished result by him and Andrew Gleason which asymptotically
improves (\ref{Elk}) with factor $\sqrt{2}$.
\end{enumerate}

\section{Sections of the cube and sinc integrals}
\label{sinc}

Let $C=[-1,1]^n\subset\R^n$ and let ${\ve
s}=(s_1,\ldots,s_n)\in\R^n$ be a unit vector. It is a well known
fact (see e. g. \cite{Ba}) that
\be \vol_{n-1}({\ve s}^\bot\cap C)
=\frac{2^n}{\pi}\int_0^\infty\prod_{i=1}^n\frac{\sin s_i t}{s_i t}\,
dt\,,\label{volume_via_sinc}\ee
where ${\ve s}^\bot$ is the $(n-1)$--dimensional subspace
orthogonal to ${\ve s}$.
In particular, the volume of the section orthogonal to the vertex
${\ve v}=(1,\ldots,1)$ of $C$ is given by
\bea \vol_{n-1}({\ve v}^\bot\cap C)
=\frac{2^n}{\pi}\int_0^\infty\left(\frac{\sin
\frac{t}{\sqrt{n}}}{\frac{t}{\sqrt{n}}}\right)^n dt=
2^{n-1}\sqrt{n}\,\sigma_n\,.\eea
Laplace and P\'{o}lya (see \cite{La}, \cite{Po} and e. g.
\cite{ChLo}) both gave proofs that
\bea \lim_{n\rightarrow\infty}\frac{\vol_{n-1}({\ve v}^\bot\cap
C)}{2^{n-1}}=\sqrt{\frac{6}{\pi}}\,. \eea
Thus, (\ref{asympt_sigma}) is justified.
\begin{lemma}
For $n\ge 2$
\bea 0<\sigma_{n+1}<\sigma_{n}\le 1\,.\eea \label{Borweins}
\end{lemma}
\begin{proof}
This result is implicit in \cite{BoBo}. Indeed, Theorem 1 (ii) of
\cite{BoBo} applied with $a_0=a_1=\ldots=a_{n}=1$ gives the
inequalities
\bea 0<\sigma_{n+1}\le\sigma_{n}\le 1\,.\eea
The strict inequality $\sigma_{n+1}<\sigma_{n}$ follows easily
from the observation that in this case the inequality in equation
(3) of \cite{BoBo} is strict with $a_{n+1} = a_0 = y =1$.
\end{proof}

\section{An application of the Busemann theorem}

Let $|\cdot|$  denote the euclidean norm. Recall that we can
associate with each 
star body $L$ the {\em distance function} 
$f_L({\ve x})=\inf\{\lambda>0 : {\ve x}\in \lambda L\}\,. $
The {\em intersection body} $IL$ of a star body $L\subset\R^n$,
$n\ge 2$ is defined as the {\ve o}--symmetric star body whose
distance function $f_{IL}$ is given by
\bea f_{IL}({\ve x})=\frac{|{\ve x}|}{\vol_{n-1}({\ve x}^\bot\cap
L)}\,.\eea
Intersection bodies played an important role in the solution to the
famous Busemann--Petty problem.
The Busemann theorem (see e. g. \cite{Ga}, Chapter 8) states that if
$L$ is {\ve o}--symmetric and convex, then $IL$ is the convex set.
This result allows us to prove the following useful inequality.
Let $f=f_{IC}$ denote the distance function of $IC$.
\begin{lemma}
For any non--zero ${\ve x}\in\R^n$
\be f\left(\frac{{\ve x}}{||{\ve x}||}\right) \le f({\ve
v})=\frac{1}{\sigma_n 2^{n-1}}\,,\label{boundary}\ee
with equality only if $n=2$ or $\frac{{\ve x}}{||{\ve x}||}$ is a
vertex of the cube $C$.
\label{distances}
\end{lemma}
\begin{proof}
%
We proceed by induction on $n$.
When $n=2$ 
the result is obvious. Suppose now (\ref{boundary}) is true for
$n-1\ge 2$. Since, if some $x_i=0$, the problem reduced to that in
$\R^{n-1}$, we may assume inductively that $x_i>0$ for all $i$.
Clearly, we may also assume that ${\ve w}=\frac{{\ve x}}{||{\ve
x}||}$ is not a vertex of $C$, in particular, ${\ve w}\neq{\ve v}$.


Let $Q=[0,1]^n\subset \R^n$ and let $L$ be the $2$--dimensional
subspace spanned by vectors ${\ve v}$ and ${\ve x}$. Then $P=L\cap
Q$ is a parallelogram on the plane $L$. To see this, observe that
the cube $Q$ is the intersection of two cones $\{{\ve y}\in\R^n:
y_i\ge 0\}$ and $\{{\ve y}\in\R^n: y_i\le 1\}$ with apexes at the
points ${\ve o}$ and {\ve v} respectively.

Suppose that $P$ has vertices ${\ve o}$, ${\ve u}$, ${\ve v}$,
${\ve v}-{\ve u}$. Then the edges ${\ve o}{\ve u}$, ${\ve o}\;{\ve
v}-{\ve u}$ of $P$ belong to coordinate hyperplanes and the edges
${\ve u}{\ve v}$, ${\ve v}\;{\ve v}-{\ve u}$ lie on the boundary
of $C$.
W. l. o. g., we may assume that the point ${\ve w}$ lies on the
edge ${\ve u}{\ve v}$.
Let
\bea {\ve v}'=\sigma_n {\ve v}=\frac{\vol_{n-1}({\ve v}^\bot\cap
C)}{2^{n-1}}\frac{{\ve v}}{|{\ve v}|}\in \frac{1}{2^{n-1}}IC\,,
\eea
\bea {\ve u}'=\sigma_{n-1}{\ve u}\,. \eea
Since the point ${\ve u}$ lies in one of the coordinate
hyperplanes, by the induction hypothesis
\bea f({\ve u}')=f(\sigma_{n-1}{\ve u})\le \frac{1}{2^{n-1}}\,.
\eea
Thus, ${\ve u}'\in \frac{1}{2^{n-1}}IC$.
Consider the triangle with vertices ${\ve o}$, ${\ve u}$, ${\ve
v}$. Let ${\ve w}'$ be the point of intersection of segments
${\ve o}{\ve w}$ and ${\ve u}'{\ve v}'$. Observing that by Lemma
\ref{Borweins}
\bea |\sigma_n{\ve w}|<|{\ve w}'|<|\sigma_{n-1}{\ve w}|\,,\eea
we get
\be \frac{1}{\sigma_{n-1}}<\frac{|{\ve w}|}{|{\ve
w}'|}<\frac{1}{\sigma_{n}}\,. \label{ratio}\ee
%
%
%
By the Busemann theorem $IC$ is convex. Therefore ${\ve
w}'\in\frac{1}{2^{n-1}}IC$ 
and thus
\bea |{\ve w}'|\le\frac{\vol_{n-1}({\ve w}^\bot\cap
C)}{2^{n-1}}\,. \eea
By (\ref{ratio}) we obtain
\bea f\left(\frac{{\ve x}}{||{\ve x}||}\right)=f({\ve w})=
\frac{|{\ve w}|}{\vol_{n-1}({\ve w}^\bot\cap C)}\le \frac{|{\ve
w}|}{2^{n-1}|{\ve w}'|}<\frac{1}{\sigma_n 2^{n-1}}\,. \eea
\end{proof}
%
%
Applying Lemma \ref{distances} to a unit vector ${\ve s}$ and
using (\ref{volume_via_sinc}) we get the following inequality for
sinc integrals.
\begin{coro}
For any unit vector ${\ve s}=(s_1,\ldots,s_n)\in\R^n$
\bea ||{\ve s}||\int_0^\infty\prod_{i=1}^n\frac{\sin s_i t}{s_i t}\,
dt\ge \int_0^\infty\left(\frac{\sin t}{t}\right)^n dt\,,\eea
with equality only if $n=2$ or $\frac{{\ve s}}{||{\ve s}||}$ is a
vertex of the cube $C$.
\end{coro}
\noindent {\bf Remark 3}
Note that $IC$ is symmetric w. r. t. any coordinate hyperplane. This
observation and Busemann's theorem immediately imply
(\ref{boundary}) with non--strict inequality in all cases.

\section{Proof of the theorem }
Clearly, we may assume that $||{\ve a}||>1$ and, in particular, that
the inequality in Lemma \ref{distances} is strict for ${\ve x}={\ve
a}$.
We shall also assume w. l. o. g. that $\gcd(a_1,\ldots,a_n)=1$.

Let $S={\ve a}^\bot\cap C$ and $\Lambda={\ve a}^\bot\cap\Z^n$. Then
$S$ is a centrally symmetric convex set and $\Lambda$ is a
$(n-1)$--dimensional sublattice of $\Z^n$ with determinant
(covolume) $\det\Lambda=|{\ve a}|$.
Let $\lambda_i=\lambda_i(S,\Lambda)$ be the $i$--th successive
minimum of $S$ w. r. t. $\Lambda$, that is
\bea \lambda_i=\inf\{\lambda> 0: \dim(\lambda S\cap\Lambda)\ge i
\}\,. \eea
By the definition of $S$ and $\Lambda$ it is enough to show that
\bea \lambda_1\cdots\lambda_{n-1}<\frac{||{\ve
a}||}{\sigma_n}\,.\eea
The $(n-1)$--dimensional subspace ${\ve a}^\bot\subset\R^n$ can be
considered as a usual $(n-1)$--dimensional euclidean space. The
Minkowski Theorem on Successive Minima (see e. g. \cite{GrLek},
Chapter 2), applied to the ${\ve o}$--symmetric convex set $S\subset
{\ve a}^\bot$ and the lattice $\Lambda\subset {\ve a}^\bot$, implies
that
\bea
\lambda_1\cdots\lambda_{n-1}\le\frac{2^{n-1}\det\Lambda}{\vol_{n-1}(S)}=\frac{2^{n-1}|{\ve
a}|}{\vol_{n-1}({\ve a}^\bot\cap C)}=2^{n-1}f({\ve a})\,, \eea
and by Lemma \ref{distances} we get
\bea \lambda_1\cdots\lambda_{n-1}\le 2^{n-1}f({\ve
a})=2^{n-1}f\left(\frac{{\ve a}}{||{\ve a}||}\right)||{\ve a}||\eea
\bea <2^{n-1}f({\ve v})||{\ve a}||=\frac{||{\ve a}||}{\sigma_n}\,.
\eea
This proves the theorem.

\section{Proof of Corollary \ref{lower}}

For a sum--distinct set $\{a_1,\ldots,a_n\}$ consider the vector
${\ve a}=(a_1,\ldots,a_n)$. Observe that any non--zero integral
vector ${\ve x}$ with $L({\ve x})=0$ must have the maximum norm
greater than 1. Therefore (\ref{product})
implies the inequality
\bea 2^{n-1}<\frac{||{\ve a}||}{\sigma_n}\,. \eea

\section{Acknowledgements}

The author wishes to thank Professors D. Borwein and A. Schinzel for
valuable comments and Professor P. Gruber for fruitful discussions
and suggestions. The work was partially supported by FWF Austrian
Science Fund, project M821--N12.


\begin{thebibliography}{99}
%
\bibitem{Disser}
I. Aliev, On a Decomposition of Integer Vectors, PhD Dissertation,
Institute of Mathematics PAN, Warsaw 2001.
%
\bibitem{Ba}
K. Ball, Cube Slicing in $\R^n$, Proc. Amer. Math. Soc. {97} (1986)
no. 3 465--472.
%
\bibitem{BombVaal}
E. Bombieri and J. Vaaler, On Siegel's Lemma, Invent. Math. {73}
(1983) 11--32, Addendum, ibid. 75 (1984) 377.
%
\bibitem{BoBo}
D. Borwein and J. Borwein, Some Remarkable Properties of Sinc and
Related Integrals, Ramanujan J. {5} (2001) no. 1 73--89.
%
\bibitem{BoMo}
P. Borwein and M. Mossinghoff, Newman polynomials with prescribed
vanishing and integer sets with distinct subset sums, Math. Comp. 72
(2003), no. 242, 787--800 (electronic).
%
\bibitem{ChLo}
D. Chakerian and D. Logothetti, Cube Slices, Pictorial Triangles,
and Probability, Math. Mag. {64} (1991) no. 4 219--241.
%
\bibitem{El}
 N. D. Elkies, An Improved Lower Bound on the Greatest Element of a Sum-Distinct Set of Fixed
 Order,
J. Combin. Theory Ser. A {41} (1986) no. 1 89--94.
%
\bibitem{EM}
P. Erd\"os, Problems and Results in Additive Number Theory, in
,,Colloque sur la Th\'eorie des Nombres, Bruxelles, 1955'', pp.
127--137.
%
\bibitem{Ga}
R. J. Gardner, Geometric Tomography, Encyclopedia of Mathematics and
its Applications, 58, Cambridge University Press, Cambridge, 1995.
%
\bibitem{GrLek}
P. M. Gruber and C. G. Lekkerkerker, Geometry of Numbers,
North--Holland, Amsterdam, 1987.
%
\bibitem{Guy}
R. K. Guy, Unsolved Problems in Number Theory, Third edition.
Problem Books in Mathematics. Unsolved Problems in Intuitive
Mathematics, Springer-Verlag, New York, 2004.
%
\bibitem{La}
P. S. Laplace, Th\'{e}orie Analytique des Probabilit\'{e}s, Paris,
1812.
%
\bibitem{MeRo}
R. G. Medhurst and J. H. Roberts, Evaluation of the Integral
$I_n(b)=\frac{2}{\pi}\int_0^\infty\left(\frac{\sin
x}{x}\right)^n\scriptstyle \cos(bx)dx$, Math. Comp. {19} (1965)
113--117.
%
\bibitem{NewSL}
A. Schinzel, A Property of Polynomials with an Application to
Siegel's Lemma, Monatsh. Math. {137} (2002) 239--251.
%
\bibitem{Po}
G. P\'{o}lya, Berechnung eines Bestimmten Integrals, Math. Ann. {
74} (1913) 204--212.
%
\bibitem{Sloane}
N. J. A. Sloane, Sequences A049330 and A049331 in {The On-Line
Encyclopedia of Integer Sequences},
http://www.research.att.com/~njas/sequences/.
%
\bibitem{Woods}
A. C. Woods, The anomaly of convex bodies. Proc. Cambridge Philos.
Soc. 52 (1956), 406--423.
%
\end{thebibliography}
\end{document}